\documentclass[11pt,letterpaper]{article}
\usepackage{amsmath}
\usepackage{amssymb}
\usepackage{verbatim}
\usepackage{graphicx}
\usepackage{theorem}

\usepackage{fullpage}

\def\B{{\mathbb B}}
\def\G{{\mathbb G}}
\def\R{{\mathbb R}}
\def\Pr{{\mathbb P}}
\def\Ex{{\mathbb E}}
\def\cB{{\mathcal B}}
\def\cF{{\mathcal F}}
\def\cM{{\mathcal M}}
\def\cP{{\mathcal P}}

\def\cW{{\mathcal W}}

\def\whp{{\bf whp}}
\def\idty{{\mathbf 1}}

\def\Erdos{Erd\H{o}s}
\def\Lovasz{Lov\'asz}

\newtheorem{theorem}{Theorem}[section]
\newtheorem{lemma}[theorem]{Lemma}

\theorembodyfont{\upshape}
\newtheorem{remark}[theorem]{Remark}
\numberwithin{equation}{section}

\newenvironment{proof}{\noindent {\sc Proof. }}{\hfill $\square$\\}

\begin{document}

\title{\bf A phase transition and stochastic domination in Pippenger's
  probabilistic failure model for Boolean networks with unreliable gates}

\author{Maxim Raginsky\thanks{Department of Electrical and Computer
     Engineering, Northwestern University, Evanston, IL 60208,
     USA. E-mail: maxim@ece.northwestern.edu}}

\date{November 24, 2003}
\maketitle

\begin{abstract}
\vskip 1em
\begin{center}
\begin{minipage}{5in}
\noindent We study Pippenger's model of Boolean networks with
unreliable gates. In this model, the conditional probability that a
particular gate fails, given the failure status of any subset of gates
preceding it in the network, is bounded from above by some
$\varepsilon$. We show that if a Boolean network with $n$
gates is selected at random according to the
Barak-\Erdos\ model of a random acyclic digraph, such that the expected edge
density is $c n^{-1}\log n$, and if $\varepsilon$ is
equal to a certain function of the size of the largest reflexive,
transitive closure of a vertex (with respect to a particular
realization of the random digraph), then Pippenger's model exhibits a
phase transition at $c=1$. Namely, with probability $1-o(1)$ as
$n\to\infty$, we have the following: for $0 \le c \le 1$, the minimum of
the probability that no gate has failed, taken over all probability
distributions of gate failures consistent with Pippenger's model, is
equal to $o(1)$, whereas for $c >1$ it is equal to
$\exp\left(-\frac{c}{e(c-1)}\right) + o(1)$. We also indicate how a
more refined analysis of Pippenger's model, e.g., for the purpose of
estimating probabilities of monotone events, can be carried out using
the machinery of stochastic domination.\\[-2em]



\end{minipage}
\end{center}
\end{abstract}

\section{Introduction}
\label{sec:intro}

The study of phase transitions in combinatorial structures, such as
random graphs \cite{AS00,Bol01,JLR00} is a subject
at the intersection of statistical physics, theoretical computer
science, and discrete mathematics
\cite{AB02,BMW00,BBCKW01,BCP01,LS97,MMZ01,PT01,SC01}. The key idea
behind this study is that large combinatorial structures can be
thought of as systems consisting of many locally interacting
components, in direct analogy to the kinds of systems within the
purview of statistical mechanics. A phase transition, then, is a
phenomenon that takes place in certain kinds of such systems in the
limit of an infinite number of components, and corresponds
qualitatively to a change in some global (macroscopic) parameter of
the system as the local parameters of the components are varied.

Boolean networks with gates subject to probabilistic failures fall
naturally into the category of systems just described. The
possibility of a phase transition arises here, for instance, when one
associates a probability of failure with each gate of such a network,
and then looks at the maximum of the probability that the network
deviates (outputs an incorrect result), taken over all possible
assignments of inputs to the network, in the limit of an infinite number of gates. The theory of
Boolean networks with unreliable gates can be traced back to the
seminal work of von Neumann \cite{Neu56}, who considered the simplest
case, namely when each gate in the network fails with fixed
probability $\varepsilon$ independently of all other gates --- we will
refer to this set-up as the {\em $\varepsilon$-independent failure model}. Von
Neumann's initial work was developed further by
Dobrushin and Ortyukov \cite{DO77a,DO77b}, Pippenger \cite{Pip88},
Feder \cite{Fed89}, Pippenger, Stamoulis, and Tsitsiklis \cite{PST91},
and G\'acs and G\'al \cite{GG94}, to name just a few. (It should be
  mentioned that in Ref.~\cite{PST91} the authors pointed out several
  technical flaws of \cite{DO77a} and presented their own proof of a
  weaker result; G\'acs and G\'al \cite{GG94} later developed methods
  to recover the full result claimed in \cite{DO77a}.) Now we will
  summarize relevant notions and ideas in a more or less narrative fashion;
  the requisite details will be supplied in Section~\ref{sec:prelims}.

One of the key results obtained by von Neumann \cite{Neu56} was the
following: if the probability $\varepsilon$ of gate failure is
sufficiently small, then any Boolean function can be computed by a
network of unreliable gates such that the probability of error is
bounded by a constant independent of the function being computed. On
the other hand, since in this model a Boolean network is no more
reliable than its last gate, the probability of error can get
arbitrarily close to one if the probability of gate failure is
sufficiently large. (This is, in fact, one possibility for a phase
transition in Boolean networks with unreliable gates, which we have alluded
to in the preceding paragraph.)

However, as pointed out by Pippenger \cite{Pip89}, the model of
independent stochastic failures has the following significant
drawback. Suppose that, within this model, a network is shown to
compute a Boolean function $f$ with probability of error at most
$\delta$, when the gate failure probability is equal to $\varepsilon$. Then it cannot
be guaranteed in general that the same network will compute $f$ with
probability of error at most $\delta$ when the gate
failure probability is smaller than $\varepsilon'$. In particular, such a network may not even compute $f$
at all in the absence of failures! This is due to the fact that gates
which fail with a fixed and known probability can be assembled into
random-number generators that would output independent and nearly
unbiased random bits. These random-number generators can, in turn,
be used to implement a randomized algorithm that would correctly
compute $f$ with high probability. However, if one were to replace the
outputs of the random-number generators with some fixed constants, then that
algorithm would be very likely to produce meaningless results. Another
observation made by Pippenger was the following. In complexity theory of Boolean circuits \cite{Weg87}, a
theorem due to Muller \cite{Mul56} says that, given any two finite,
complete bases $\cB$ and $\cB'$ of Boolean functions, a network over
$\cB$ that computes a function $f$ can be realized as a network over
$\cB'$ with size and depth differing from those of the original
circuit by multiplicative constants that depend only on $\cB$ and
$\cB'$. It is not immediately clear under what conditions such an invariance
theorem would hold for networks with unreliable gates.

In order to overcome these objections, Pippenger proposed in
\cite{Pip89} a more general model of Boolean networks with unreliable
gates. Gate failures under this model are no longer independent, but instead are
such that the conditional probability of any gate failing, given the
status (failed or not) of any set of gates preceding it, is at most
$\varepsilon$. In Pippenger's terminology \cite{Pip89}, this model is
called {\em $\varepsilon$-admissible}. It is immediately evident that
the $\varepsilon$-admissible model subsumes the
$\varepsilon$-independent one. It also follows from definitions that a
network that computes a function $f$ reliably for all probability
distributions of gate failures within the
$\varepsilon$-admissible model, will continue to do so under the
$\varepsilon'$-admissible model for any $0 \le \varepsilon' \le
\varepsilon$. Another key achievement of Pippenger's paper \cite{Pip89} is
the proof of a Muller-type invariance theorem for Boolean networks
with unreliable gates, in
which the $\varepsilon$-admissible model plays an essential role.

The contribution of the present paper consists mainly in showing that
Pippenger's $\varepsilon$-admissible model, applied to Boolean
networks drawn at random according to a certain model of random
directed acyclic graphs, exhibits a phase transition in terms of the
minimum probability of the failure-free configuration as the
network's wiring pattern evolves from ``sparse'' to ``dense.'' The paper is organized as follows. In Section~\ref{sec:prelims} we
fix definitions and notation used throughout the paper and collect
some preliminaries on graphs, Boolean networks, and the formalism used
in Pippenger's paper \cite{Pip89}. Our main result --- one concerning the
phase transition --- is proved in Section~\ref{sec:phasetrans}. Then,
in Section~\ref{sec:stochdom}, we use the machinery of stochastic
domination to carry out a systematic analysis of the more delicate
features of Pippenger's model. We close with some remarks in
Section~\ref{sec:close} concerning directions for future
research. Finally, in the Appendix we prove a certain theorem which,
though somewhat tangential to the matter at hand, is closely related to
some mathematical techniques and concepts used in this paper.

\section{Preliminaries, definitions, notation}
\label{sec:prelims}

\subsection{Graphs}
\label{ssec:graphs}

In this paper we deal exclusively with directed acyclic graphs
(or acyclic digraphs). Given
such a graph $G=(V,E)$, we will follow standard practice of denoting
by $v(G)$ the number of vertices of $G$ and by $e(G)$ the number of
edges. Any acyclic digraph has at least one vertex of in-degree
zero. We will denote by $\hat{G}=(\hat{V},\hat{E})$ the graph obtained
from $G$ by deleting all such vertices along with all of their outgoing edges.

Let us define the {\em out-neighborhood} of a vertex $i \in V$ as the
set $N(i) := \{ j \in V: (i,j) \in E\}$, and the {\em closed
  out-neighborhood} as $\bar{N}(i) := N(i) \cup \{ i\}$. If vertex $j$ can be reached from vertex $i$ by a directed path, we
will write $i \rightsquigarrow j$ (or $i \rightsquigarrow_G j$
whenever we need to specify $G$ explicitly). The {\em transitive
  closure} of a graph $G$ is the graph $G^\star = (V^\star,E^\star)$
with $V^\star = V$ and $E^\star = \{ (i,j) : i
\rightsquigarrow_G j \}$. The transitive closure of a vertex $i$ is
the set $\Gamma(i) = \{ j \in V : (i,j) \in E^\star \}$; the set
$\Gamma^\star(i) = \Gamma(i) \cup \{ i\}$ is called the {\em
  reflexive, transitive closure} (RTC, for short) of $i$. Note that $\Gamma^{(\star)}(i)$ is
simply the (closed) out-neighborhood of the vertex $i$ in $G^\star$. It is also convenient to partially order the
vertices of $G$ as follows: for $i,j$ distinct, we will write $i
\preceq j$ if $i \rightsquigarrow j$, and require that $i \preceq i$
for each $i$. In this way, $\preceq$ is simply the reflexive closure
of the asymmetric transitive relation $\rightsquigarrow$. 

An important role in this paper will be played by the random acyclic digraph introduced by Barak and \Erdos\
\cite{BE84}. It is obtained from the standard undirected random graph
$\G(n,p)$ \cite{AS00,Bol01,JLR00} by orienting the edges according to the natural ordering of the
vertex set $[n]$, and will henceforth be denoted by $\G_d(n,p)$. 

\subsection{Boolean networks}
\label{ssec:boole}

A {\em Boolean function} is
any function $f : \B^s \rightarrow \B$, where $\B := \{ 0,1\}$. A set
of Boolean functions is referred to as a {\em basis}. In particular, we
say that a basis $\cB$ is {\em complete} if any Boolean function can
be realized by composing elements of $\cB$. Let $\cB$ be a finite
complete basis. A {\em Boolean network} (or {\em circuit}) ${\sf N}$ {\em over
  $\cB$} is an acyclic digraph $G$ with a specially designated
vertex of out-degree zero (the {\em output} of ${\sf N}$), such that each vertex of
$\hat{G}$ is labelled by some Boolean function $\varphi \in
\cB$ of its immediate predecessors, and each vertex in $V \backslash \hat{V}$
is labelled either by a Boolean variable (these vertices are the {\em
  inputs} of ${\sf N}$) or by a constant 0 or 1. Whenever there is a
need to specify the network ${\sf N}$ explicitly, we will write, e.g.,
$G_{\sf N} = (V_{\sf N},E_{\sf N})$, etc. We will refer to the graph
$\hat{G}_{\sf N}$ as the {\em gate interconnection graph} of ${\sf
      N}$. Given a network ${\sf N}$ with $s$
input vertices, we will assume the latter to be ordered in some way,
and therefore $x_i$, $1 \le i \le s$, will denote the Boolean variable
associated with the $i$th input vertex. For any assignment
$(x_1,\ldots,x_s) \in \B^s$ of values to the inputs of the network,
the value of each vertex can be computed recursively in the obvious
way, namely by evaluating the Boolean function labelling it on the
values of its immediate predecessors. We then say that the network
{\em computes} a Boolean function $f : \B^s \rightarrow \B$ if, for
any $(x_1,\ldots,x_s) \in \B^s$, the value of the output vertex, which
we will denote by ${\sf N}(x_1,\ldots,x_s)$, is equal to
$f(x_1,\ldots,x_s)$.

Let us associate with a network ${\sf N}$ the measurable space
$(\Omega_{\sf N},\cF_{\sf N})$, where $\Omega_{\sf N} :=
\B^{\hat{V}_{\sf N}}$ and $\cF_{\sf N}$ is the set of all subsets of
$\Omega_{\sf N}$. Then the occurrence of failures in the gates of
${\sf N}$ is described by a probability measure $\mu$ on $(\Omega_{\sf
  N},\cF_{\sf N})$ or, equivalently, by a family $\{X_i : i \in
  \hat{V}_{\sf N}\}$ of $\B$-valued random variables, where $X_i$ is
  the indicator function of the event
\begin{equation}
A_i := \{ \text{gate $i$ fails}
  \} \equiv \{ x \in \B^{\hat{V}_{\sf N}} : x_i =1 \},
\end{equation}
and the equivalence is, of course, given by
\begin{equation}
\Pr \big(\bigwedge_{i \in Y}(X_i = 1) \big) = \mu \big( \bigcap_{i \in
  Y}A_i \big) \qquad Y \in \cF_{\sf N}.
\end{equation}
From now on, given a probability measure $\mu$, we will denote
  probabilities of various events by $\mu(\cdot)$ or by $\Pr_\mu(\cdot)$,
  or sometimes by just $\Pr(\cdot)$, whenever the omission of the
  underlying measure is not likely to cause ambiguity.

Following
  Pippenger \cite{Pip89}, we define a {\em probabilistic failure
    model} (or PFM, for short) as a map $M$ that assigns to every Boolean network ${\sf
    N}$ a compact subset $M({\sf N})$ of the set $\cP({\sf N})$ of
  all probability measures on $(\Omega_{\sf N},\cF_{\sf N})$. One
  typically works with a family \{$M_\varepsilon :0 \le \varepsilon
  \le 1\}$ of PFM's, where $\varepsilon$ can be thought of as
  a local parameter describing the behavior of individual gates; to
  give a simple example, the $\varepsilon$-independent PFM is the map $M_\varepsilon$ that assigns to each
  network ${\sf N}$ the product measure $\pi^{\sf N}_\varepsilon :=
  \prod_{i \in \hat{V}_{\sf  N}}\nu^i_\varepsilon$, where each $\nu^i_\varepsilon$ is a copy
  of the Bernoulli measure $\nu$ with $\nu(1)=\varepsilon$. Given such a family ${\bf M} := \{ M_\varepsilon \}$, a network ${\sf N}$ with $s$ inputs and a Boolean function $f : \B^s
  \to \B$, we say that ${\sf N}$ $(\varepsilon,\delta)$-{\em
    computes $f$ with respect to ${\bf M}$} if
\begin{equation}
\max_{(x_1,\ldots,x_s) \in \B^s} \sup_{\mu \in
  M_\varepsilon({\sf N})}\Pr \big({\sf N}(x_1,\ldots,x_s) \neq
  f(x_1,\ldots,x_s)\big) \le \delta.
\label{eq:perror}
\end{equation}
The maximum in Eq.~(\ref{eq:perror}) exists owing to the finiteness of
$\B^s$ and to the compactness of $M_\varepsilon({\sf N})$. Whenever
the family ${\bf M}$ contains only one PFM $M$, we will assume that
the underlying parameter $\varepsilon$ is known and fixed, and say
that ${\sf N}$ $(\varepsilon,\delta)$-computes $f$ with respect to $M$.

Consider a pair of PFM's, $M$ and $M'$. In the terminology of Pippenger
\cite{Pip89}, $M$ is {\em more stringent} than $M'$ if, for any
network ${\sf N}$, $M({\sf N}) \supseteq M'({\sf N})$. We will also
say that $M'$ is less stringent than $M$. Thus, if one is
able to show that a network ${\sf N}$ $(\varepsilon,\delta)$-computes a
function $f$ with respect to a PFM $M$, then the same network will
also $(\varepsilon,\delta)$-compute $f$ with respect to any PFM $M'$ less
stringent than $M$.

We would also like to comment on an interesting ``adversarial''
aspect of the PFM formalism (see also Ref.~\cite{Pip90}). Let us fix a
family $\{ M_\varepsilon \}$ of PFM's. We can then envision the following game
played by two players, the Programmer and the Hacker, with the aid of
a disinterested third party, the Referee. The Referee picks a constant
$\varepsilon_0 \in (0,1)$ and announces it to the players. The Programmer picks a Boolean function $f$ and designs a
network ${\sf N}$ that would compute $f$ in the absence of
failures. He then presents ${\sf N}$ to the Hacker and lets him choose (a)
the input to ${\sf N}$ and (b) the locations of gate failures
according to $M_{\varepsilon_0}$. We assume here that the Hacker
possesses full knowledge of the structure of ${\sf N}$. The Hacker's
objective is to force the network to $(\varepsilon_0,\delta)$-compute
$f$ with $\delta > 1/2$, and the Programmer's objective is to design
${\sf N}$ in such a way that it $(\varepsilon_0,\delta')$-computes $f$ with
$\delta' < 1/2$, regardless of what the Hacker may do.

\subsection{Pippenger's model}
\label{ssec:pippenger}
Now we state the precise definition of the $\varepsilon$-admissible
PFM of Pippenger \cite{Pip89}, alluded to in
the Introduction. Given a network ${\sf N}$, let $\cM_\varepsilon({\sf
  N})$ be the set of all probability measures $\mu \in \cP({\sf N})$
that satisfy the following condition: for any gate $i \in \hat{V}_{\sf
  N}$ and for any two disjoint sets $Y,Y' \subseteq \hat{V}_{\sf N} \backslash
\Gamma^\star(i)$, such that
\begin{equation}
\mu\big(\bigcap_{j \in Y}A_j \cap \bigcap_{j \in Y'}\overline{A_j}
\big) \neq 0,
\end{equation}
we have
\begin{equation}
\mu \big( A_i | \bigcap_{j \in Y}A_j \cap \bigcap_{j \in Y'}
\overline{A_j} \big) \le \varepsilon.
\label{eq:admiss}
\end{equation}
According to definitions set forth in Section~\ref{ssec:boole}, we will have a PFM ${\sf
  N} \mapsto \cM_\varepsilon({\sf N})$ if we prove that
  $\cM_\varepsilon({\sf N})$ is a compact set. This is accomplished in
  the lemma below (incidentally, this issue has not been
  addressed in Pippenger's paper \cite{Pip89}).

\begin{lemma}\label{lm:admcompact} The set $\cM_\varepsilon({\sf N})$
  is compact in the metric topology induced by the total variation
  distance \cite{Dur96}
\begin{equation}
d(\mu,\mu') := \sup_{A \in \cF_{\sf N}} |\mu(A) - \mu'(A)|.
\end{equation}
\end{lemma}

\begin{remark}The topology induced by the total variation distance is
  actually a norm topology, where the total variation norm is defined
  on the set $\cM_\pm({\sf N})$ of all {\em signed} measures on $(\Omega_{\sf
  N},\cF_{\sf N})$ by
\begin{equation}
\|\mu\| := \sup_{A \in \cF_{\sf N}}|\mu(A)|.
\end{equation}
Furthermore, $\cP({\sf N})$, obviously being closed and bounded with respect to
the total variation norm, is a compact subset of
$\R^{v(\hat{G}_{\sf N})}$. 
\end{remark}

\begin{proof} Since $\cP({\sf N})$ is compact (see Remark above), it suffices to show
  that $\cM_\varepsilon({\sf N})$ is closed. Suppose that a sequence $\{ \mu_n \}$ in
  $\cM_\varepsilon({\sf N})$ converges to $\mu$ in total variation
  distance. Let us adopt the following shorthand
  notation: for any two disjoint sets $Y,Y' \subseteq \hat{V}_{\sf N}$, let
\begin{equation}
C_{Y,Y'} := \bigcap_{j \in Y}A_j \cap \bigcap_{j \in
  Y'}\overline{A_j}.
\end{equation}
Fix a gate $i$. Let $Y,Y' \subseteq \hat{V}_{\sf N} \backslash
\Gamma^\star(i)$ be disjoint sets such that $\mu(C_{Y,Y'}) \neq
0$. Then we can find a subsequence $\{ \mu_{n_\alpha} \}$, such that each
$\mu_{n_\alpha}(C_{Y,Y'})$ is nonzero as well. By
$\varepsilon$-admissibility, we have the following estimate:
\begin{eqnarray}
\mu (A_i | C_{Y,Y'} ) &\le& |\mu(A_i | C_{Y,Y'}) -
\mu_{n_\alpha}(A_i | C_{Y,Y'})| + \mu_{n_\alpha}(A_i | C_{Y,Y'}) \nonumber
\\
&\le&  |\mu(A_i | C_{Y,Y'}) -
\mu_{n_\alpha}(A_i | C_{Y,Y'})| + \varepsilon.\label{eq:compstep1}
\end{eqnarray}
We can further estimate the first term on the right-hand side of
(\ref{eq:compstep1}):
\begin{eqnarray}
&& |\mu(A_i | C_{Y,Y'}) -
\mu_{n_\alpha}(A_i | C_{Y,Y'})| = \left| \frac{\mu(A_i \cap C_{Y,Y'})}{\mu(C_{Y,Y'})} -
\frac{\mu_{n_\alpha}(A_i \cap
  C_{Y,Y'})}{\mu_{n_\alpha}(C_{Y,Y'})}\right| \nonumber \\
&& \qquad \le \frac{1}{\mu(C_{Y,Y'})}\big|\mu(A_i \cap C_{Y,Y'}) -
\mu_{n_\alpha}(A_i \cap C_{Y,Y'})\big| \nonumber \\
&& \qquad \qquad + \frac{\mu_{n_\alpha}(A_i |
  C_{Y,Y'})}{\mu(C_{Y,Y'})} \left|\mu_{n_\alpha}(C_{Y,Y'}) - \mu(C_{Y,Y'})\right|
\nonumber \\
&& \qquad \le \frac{1+\varepsilon}{\mu(C_{Y,Y'})}d(\mu,\mu_{n_\alpha}).\label{eq:compstep2}
\end{eqnarray}
Combining (\ref{eq:compstep1}) and (\ref{eq:compstep2}) and taking the
limit along $n_\alpha$, we obtain $\mu(A_i | C_{Y,Y'} ) \le
\varepsilon$. Thus $\cM_\varepsilon({\sf N})$ is closed, hence compact.
\end{proof}

\noindent As we have mentioned earlier, Pippenger \cite{Pip89} has termed the PFM
${\sf N} \mapsto \cM_\varepsilon({\sf N})$ {\em
  $\varepsilon$-admissible}. We will also abuse language slightly by
referring to individual probability measures $\mu \in \cM_\varepsilon({\sf
  N})$ as $\varepsilon$-admissible.

It is easy to see that $\cM_\varepsilon({\sf N})$ containts all
Bernoulli product measures $\pi^{\sf N}_{\varepsilon'}$ with $0 \le
\varepsilon' \le \varepsilon$, as well as all product measures
$\prod_{i \in \hat{V}_{\sf N}}\nu^i_{\varepsilon_i}$ with $0 \le
\varepsilon_i \le \varepsilon$. Furthermore, it follows directly from
definitions that $\cM_{\varepsilon'}({\sf
  N}) \subseteq \cM_\varepsilon({\sf N})$ for $0 \le \varepsilon' \le
\varepsilon$. That is, the $\varepsilon$-admissible PFM is more stringent than the
$\varepsilon'$-admissible one. Therefore, when $\varepsilon' \in
     [0,\varepsilon]$, a network that $(\varepsilon,\delta)$-computes
a function $f$ under the $\varepsilon$-admissible model will also
$(\varepsilon',\delta)$-compute the same function under the
$\varepsilon'$-admissible model and, in particular, when the gate
failures are distributed according to $\pi^{\sf N}_{\varepsilon'}$.

\section{The phase transition}
\label{sec:phasetrans}

\subsection{Motivation and heuristics}
\label{ssec:motivation}

Our main result, to be stated and proved in the next section, is
formulated in terms of the probability of the failure-free configuration in a
network of unreliable gates, under the $\varepsilon$-admissible model
of Pippenger. As we shall demonstrate shortly, this quantity does not
depend on the particular function being computed, but only on the size
and the structure of the gate interconnection graph associated to the
network.

Given a Boolean network ${\sf N}$, let us consider the quantity
\begin{equation}
\inf_{\mu \in \cM_{\varepsilon}({\sf N})} \mu \big(\bigcap_{i \in
  \hat{V}_{\sf N}}\overline{A_i} \big).
\label{eq:ff}
\end{equation}
The set $\cM_\varepsilon({\sf N})$ is compact by
Lemma~\ref{lm:admcompact}, and the expression being minimized is a
continuous function of $\mu$ with respect to total variation
distance. Thus, the infimum in (\ref{eq:ff}) is actually attained, and
a moment of thought reveals that this quantity depends
only on the structure of the gate interconnection graph of ${\sf N}$,
but not on the specific gate labels or on the identity of the output
vertex. Therefore, given an acyclic digraph $G=(V,E)$, let us define
$F_\varepsilon(G)$ as the quantity (\ref{eq:ff}) for all networks
${\sf N}$ whose gate interconnection graphs are isomorphic to $G$, modulo
gate labels and the identity of the output vertex. In the same spirit,
let us denote by $\cM^G_\varepsilon$ the set of all
$\varepsilon$-admissible probability measures on the measurable space
$(\Omega_G,\cF_G)$ where, as before, $\Omega_G := \B^V$ and $\cF_G$ is
the set of all subsets of $\Omega_G$. Then we can write
\begin{equation}
F_\varepsilon(G) := \inf_{\mu \in \cM^G_\varepsilon} \mu \big(
\bigcap_{i \in V}\overline{A_i} \big).
\label{eq:ff1}
\end{equation}
Our motivation to focus on $F_\varepsilon(G)$ is twofold: firstly, we
are able to gloss over such details as the function being
computed or the basis of Boolean functions used to construct a given
network, and secondly, $F_\varepsilon(G)$ can also be used to obtain
lower bounds on probabilities of other events one would associate with
``proper'' operation of the network (such as, e.g., the probability that the majority of gates have not failed \cite{Pip89}).

In order to get a quick idea of the dependence of $F_\varepsilon(G)$ on the
structure of $G$, we can appeal to the \Lovasz\ local
lemma \cite{EL75} or, rather, to a variant thereof due to \Erdos\ and
Spencer \cite{ES91}. (See also Alon and Spencer \cite{AS00} and
Bollob\'as \cite{Bol01} for proofs
and a sampling of applications.) The basic idea behind the \Lovasz\
local lemma is this: we have a finite family $\{ H_i \}$ of ``bad''
events in a common probability space, and we are interested in the
probability that none of these events occur, i.e., $\Pr \big(\bigcap_i
\overline{H_i}\big)$. The ``original'' local lemma \cite{EL75} gives a
sufficient condition for this probability to be strictly positive when
most of the events $H_i$ are independent, but with strong
dependence allowed between some of the subsets of $\{ H_i
\}$; for this reason it is formulated in terms of the dependency
digraph \cite{AS00,Bol01} for $\{H_i\}$. The version due to \Erdos\
and Spencer \cite{ES91} (often referred to as ``lopsided \Lovasz\
local lemma'') has the same content, but under the weaker condition that certain
conditional probabilities involving the $H_i$ and their complements
are suitably bounded. More precisely:

\begin{theorem}[\Erdos\ and Spencer \cite{ES91}]\label{thm:4l} Let $\{ H_i \}^n_{i=1}$ be a
  family of events in a common probability space. Suppose that there
  exist a directed graph $G=(V,E)$ with $v(G)=n$ and real
  constants $\{ r_i \}^n_{i=1}$, $0 \le r_i < 1$, such that, for any
  $Y \subseteq V\backslash \bar{N}(i)$,
\begin{equation}
\Pr\big(H_i | \bigcap_{j \in Y}\overline{H_j} \big) \le r_i \prod_{j
  \in \bar{N}(i)}(1-r_j).\label{eq:4l1}
\end{equation}
Then
\begin{equation}
\Pr\big( \bigcap^n_{i=1}\overline{H_i} \big) \ge \prod^n_{i=1}(1-r_i)
> 0.\label{eq:4l2}
\end{equation}
In other words, the event ``none of the events $H_i$ occur'' holds
with strictly positive probability.
\end{theorem}

Consider now an acyclic digraph $G=(V,E)$. Then, by defintion of
$\varepsilon$-admissibility, every probability measure $\mu \in
\cM^G_\varepsilon$ satisfies
\begin{equation}
\mu\big( A_i | \bigcap_{j \in Y}\overline{A_j}\big) \le \varepsilon,
  \qquad \forall Y \subseteq V\backslash
  \Gamma^\star(i).
\label{eq:pipcond}
\end{equation}
We can rewrite (\ref{eq:pipcond}) in terms of the transitive closure
graph $G^\star$ as follows. Denote the out-neighborhood of a vertex
$i$ in $G^\star$ by $N^\star(i)$, and similarly for the closed
out-neighborhood. Then (\ref{eq:pipcond}) becomes
\begin{equation}
\mu\big( A_i | \bigcap_{j \in Y}\overline{A_j}
\big) \le \varepsilon, \qquad \forall Y \subseteq V^\star \backslash \bar{N}^\star(i) .
\end{equation}
Now let $\Delta$ be the maximum out-degree of $G^\star$, i.e., $\Delta
:= \max_{i \in V} |\Gamma(i)|$. Then, provided that $\varepsilon \le
\Delta^\Delta/(\Delta+1)^{\Delta + 1}$, the events $\{ A_i : i \in
V\}$ will satisfy the condition (\ref{eq:4l1})
of Theorem~\ref{thm:4l} with $r_i = 1/(\Delta + 1)$ for all $i$, for
every $\mu \in \cM^G_\varepsilon$. Using (\ref{eq:4l2}), we conclude that
\begin{equation}
F_\varepsilon(G) \ge \left(1-\frac{1}{\Delta+1}\right)^{v(G)}, \qquad
\varepsilon \le \Delta^\Delta/(\Delta+1)^{\Delta+1}.
\label{eq:fgbound0}
\end{equation}
Furthermore, $F_{\varepsilon'}(G) \ge F_\varepsilon(G)$ for
$\varepsilon' \le \varepsilon$ because then we have
$\cM^G_{\varepsilon'} \subseteq \cM^G_\varepsilon$. Thus, defining $F_G :=
F_{\Delta^\Delta/(\Delta + 1)^{\Delta + 1}}(G)$, we can write
\begin{equation}
F_\varepsilon(G) \ge F_G \ge \left(1-\frac{1}{\Delta+1}\right)^{v(G)}, \qquad
\varepsilon \le \Delta^\Delta/(\Delta+1)^{\Delta+1}.
\label{eq:fgbound}
\end{equation}
In fact, Lemma~\ref{lm:f_exact} in the next section can be used to
obtain the exact expression
\begin{equation}
F_G \equiv \left(1 -
\frac{\Delta^\Delta}{(\Delta+1)^{\Delta+1}}\right)^{v(G)},
\label{eq:fgexact}
\end{equation}
i.e., $F_G$ is equal precisely to the probability of $\bigcap^{v(G)}_{i=1}\overline{A_i}$ when
the $A_i$ are independent and $\Pr(A_i) = \Delta^\Delta/(\Delta +
1)^{\Delta+1}$.

An inspection of the form of (\ref{eq:fgexact}) suggests the
following strategy for exhibiting a phase transition: We will consider a
suitable parametrization $p_c(n)$ of the (average)
density\footnote{The {\em density} of a graph $G$ is defined as $d(G)
  := e(G)/v(G)$. Thus, we have for the
expected density of the random graph $\G(n,p)$
$$
\Ex d(\G(n,p)) = \sum^n_{k=0}{n \choose k}\frac{k}{n}p^k (1-p)^{n-k} = p.
$$
The expected density of the random acyclic digraph $\G_d(n,p)$
is the same.} of the random graph $\G_d(n,p_c(n))$, with $p_c(n) \to 0$ as $n\to \infty$, such that the graph is ``sparse'' for $c <
1$ and ``dense'' for $c > 1$. Furthermore, this change from ``sparse''
to ``dense'' will be accompanied by a phase transition manifesting itself in
the distinct large-$n$ behavior of the size $\gamma^\star_n$ of the largest RTC of a vertex in $\G_d(n,p_c(n))$ depending on
whether $c < 1$, $c=1$, or $c > 1$, respectively. (This phase
transition was discovered and studied by Pittel and Tungol
\cite{PT01}, and will play a key role in our proof.) Given a
particular realization of $\G_d(n,p_c(n))$, $\gamma^\star_n = \Delta +
1$. Defining random variables $\varepsilon_n :=
\vartheta(\gamma^\star_n)$ [we have defined $\vartheta(x) :=
(x-1)^{x-1}/x^x$ in order to avoid cluttered equations], we
will end up with a sequence of $\varepsilon_n$-admissible PFM's, such
that the asymptotic behavior of $F_{\varepsilon_n}(\G_d(n,p_c(n))$
will be different depending on whether $c$ is above or below unity.

Of course, the class of probability measures satisfying the
condition (\ref{eq:pipcond}) is much broader than $\cM^G_\varepsilon$.
In terms of Boolean networks, it describes the probabilistic
failure model under which the conditional probability for a particular
gate to fail, given that (any subset of) the gates preceding it have {\em
  not} failed, is at most $\varepsilon$. (One can use the strategy of
Lemma~\ref{lm:admcompact} to prove that the
corresponding set of probability measures is compact.) In order to get a better grip
on the $\varepsilon$-admissible model, one has to make full use of its
definition; this will be done in Section~\ref{sec:stochdom} using the
machinery of stochastic domination \cite{Lig85,LSS97,Lin92}. As far as the
results in the next section are concerned, though, the condition
(\ref{eq:pipcond}) is all that is needed.\footnote{Ironically enough,
if at the very outset we were to use stochastic domination formalism to analyze $\mu\big(\bigcap_{i \in V(G)}\overline{A_i}\big)$, $\mu \in \cM^G_\varepsilon$, we would not have been able to spot the
role of the RTC in our development.} Incidentally, it is possible
to prove a more specialized version of the lopsided \Lovasz\ local
lemma which, among other things, gives a sufficient condition to have $\Pr
\big(\bigcap_i \overline{E_i}\big) > 0$ when the conditional
probabilities
\begin{equation}
\Pr \big(H_i | \bigcap_{j \in Y}H_j \cap \bigcap_{j \in
  Y'} \overline{H_j}\big)
\end{equation}
for all disjoint $Y,Y' \subseteq V\backslash \bar{N}(i)$ are suitably
bounded. Since this result is, strictly speaking, tangential to the
matter of this paper, we return to it in the Appendix.

\subsection{The main result}
\label{ssec:mainresult}

\begin{figure}
\centerline{\includegraphics[width=0.5\textwidth]{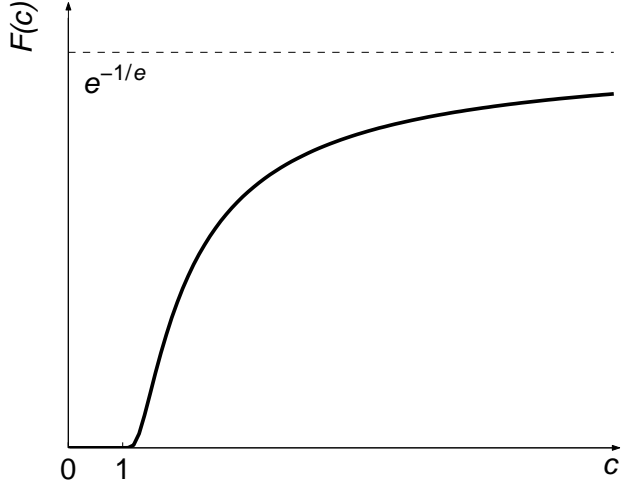}}
\caption{The phase transition of Theorem \ref{thm:main}: in
  the subcritical phase $(0 \le c < 1)$, the minimum probability of
  failure-free operation is zero; in the supercritical phase $(c > 1)$, it is
  $\exp\left(-\frac{c}{e(c-1)}\right)$. Note that $F(c)$ approaches
  $e^{-1/e}$ as $c\to \infty$.}
\label{fig:phasetrans}
\end{figure}

Now we are in a position to state and prove the main result of this
paper (the notation we use has been defined in the preceding section):

\begin{theorem}\label{thm:main} Consider all Boolean networks with $n$
  gates whose gate interconnection graphs are drawn from the random
  acyclic digraph $\G_d(n,c n^{-1}\log n)$. Let
  $\gamma^\star_n(c)$ denote the
  size of the largest RTC of a vertex in
  $\G_d(n,cn^{-1}\log n)$. Define the sequence of random
  variables
\begin{equation}
F_n(c) := F_{\vartheta(\gamma^\star_n(c))}\big(\G_d(n,c
  n^{-1}\log n)\big).
\end{equation}
Then \whp\footnote{We follow
  standard practice of writing that a sequence of events $E_n$ occurs
  \whp\ (with high probability) if $\Pr(E_n) \to 1$ as $n \to \infty$.} 
\begin{equation}
F_n(c) \rightarrow F(c) \equiv \begin{cases}
0 & \text{if $c \le 1$}\\
\exp\left(-\frac{c}{e(c-1)}\right) & \text{if $c > 1$}
\end{cases}.\label{eq:phasetrans}
\end{equation}
The corresponding phase transition is illustrated in
Fig.~\ref{fig:phasetrans}.
\end{theorem}

\begin{remark}The proof of the theorem goes through [modulo obvious
  modifications involving the exponent in the $c > 1$ case of
  (\ref{eq:phasetrans})] if, instead of $\vartheta(x)$, we use any
  nonnegative function $f(x)$ that behaves like $1/x$ for large $x$. 
\end{remark}

\begin{proof} In order to carry on, we first need to gather some
  preliminary results. Once we have all the right pieces in place, the proof is
  actually surprisingly simple.

The first result we need is the following exact formula for
$F_\varepsilon(G)$:

\begin{lemma}\label{lm:f_exact} For any acyclic digraph $G=(V,E)$ with $v(G) =
  n$, $F_\varepsilon(G) = (1-\varepsilon)^n$.
\end{lemma}

\begin{proof} Let $i_1 < i_2 < \ldots < i_n$ be an arrangement of the
  vertices of $G$ according to some linear extension \cite{Aig79} of the partial
  order $\preceq$ defined in Section~\ref{ssec:graphs}. Then, for each $k
  \in [n]$, 
\begin{equation}
i_j \not\succeq i_k, \qquad 1 \le j < k,
\end{equation}
so, by definition of
  $\varepsilon$-admissibility, we have for any $\mu \in \cM^G_\varepsilon$
\begin{equation}
\mu\big(A_{i_{k+1}} | \bigcap\nolimits^k_{j=1}\overline{A_{i_j}} \big) \le \varepsilon,
\qquad 0 \le k \le n-1.
\end{equation}
Then
\begin{equation}
\mu\big( \bigcap_{i \in V} \overline{A_i} \big) = \big(1-\mu(A_{i_1})\big)\times \prod^{n-1}_{k=1}
\left(1-\mu\big(A_{i_{k+1}} | \bigcap\nolimits^k_{j=1} \overline{A_{i_j}}\big)\right) \ge (1-\varepsilon)^n.
\label{eq:ffbound}
\end{equation}
The choice $\mu = \pi^V_\varepsilon \equiv \prod_{i \in V}\nu^i_\varepsilon$ attains the bound in
(\ref{eq:ffbound}), and the lemma is proved.
\end{proof}

Next, we will need a result on the size of the largest reflexive,
transitive closure of a vertex in the random acyclic digraph
$\G_d(n,c n^{-1} \log n)$. Pittel and Tungol \cite{PT01} showed that
the following phase transition takes place as $c$ is varied:

\begin{lemma}\label{lm:reftrancl} For the
  random acyclic digraph $\G_d(n,cn^{-1}\log n)$ one has
  the following:
\begin{enumerate}
\item If $c \ge 1$, then there exists a positive constant $A(c)$, such
  that
\begin{equation}\label{eq:pt1}
\lim_{n \to \infty} \Pr \left\{ \left| \gamma^\star_n - n
  \left(1-\frac{1}{c}\right) - \frac{2 n \log \log n}{c \log n}
  \right| \le \frac{A(c)n}{\log n}\right\} = 1.
\end{equation}
\item If $c < 1$, then for every $\varkappa > 0$, 
\begin{equation}\label{eq:pt2}
\lim_{n \to \infty} \Pr \left\{ (1-\varkappa)n^c \log n <
  \gamma^\star_n \le n^c \log n \right\} = 1.
\end{equation}
\end{enumerate}
\end{lemma}

\noindent{The exact (and rather unwieldy) expressions of Lemma
  \ref{lm:reftrancl} are more than is needed for the purposes of the
  present proof. We will settle for a rather more prosaic asymptotic form that
  follows directly from (\ref{eq:pt1}) and (\ref{eq:pt2}). Namely, \whp
\begin{equation}
\gamma^\star_n(c) \sim \begin{cases}
n^c \log n & \text{if $0 \le c < 1$} \\
\\[-1em]
\frac{2 n \log \log n}{\log n} & \text{if $c = 1$}\\
\\[-1em]
n\left(1-\frac{1}{c}\right) & \text{if $c > 1$}
\end{cases},\label{eq:reftrancl}
\end{equation}
where, as usual, the notation $a_n \sim b_n$ means that $a_n = b_n
\big(1+o(1)\big)$.}

Next, we need asymptotics of the function
$\vartheta(x)$ for large $x$. To that end, we write
\begin{equation}
\lim_{x \to\infty} x \vartheta(x) = \lim_{x \to \infty} \left(1-\frac{1}{x+1}\right)^x = \frac{1}{e}.\label{eq:asymtheta}
\end{equation}
In other words, $\vartheta(x) \sim \frac{1}{ex}$ as $x \to
\infty$.

Finally, we will have use for the following two limits:

\begin{lemma}\label{lm:limits}
\begin{eqnarray}
\lim_{n \to \infty}\left(1-\frac{1}{n^c \log n}\right)^n &=& 0, \qquad
0\le c <1\label{eq:limit1} \\
\lim_{n \to \infty}\left(1-\frac{\log n}{n \log\log n}\right)^n &=& 0.\label{eq:limit2}
\end{eqnarray}
\end{lemma}

\begin{proof} We will prove (\ref{eq:limit1}); the same strategy will
  work for (\ref{eq:limit2}) as well. Let us assume that $\log$ denotes
  natural logarithms [otherwise we can multiply the second term in
  parentheses of (\ref{eq:limit1}) by a suitable constant]. It
  suffices to show that
\begin{equation}
\lim_{n \to \infty}n \log\left(1-\frac{1}{n^c \log n}\right) =
-\infty.
\label{eq:neginfty}
\end{equation}
Using the inequality $\log x \le x-1$, we have for each $n$
\begin{equation}
n \log\left(1-\frac{1}{n^c \log n}\right) \le - \frac{n^{1-c}}{\log
  n}.
\label{eq:limit_estimate}
\end{equation}
Given any $K > 0$, one can find large enough $N$ such that the
right-hand side of (\ref{eq:limit_estimate}) is less than $-K$ for all
$n \ge N$. Therefore the same holds for the left-hand side, and
(\ref{eq:neginfty}) is proved.
\end{proof}

The rest is fairly straightforward. Using Lemma \ref{lm:f_exact} and (\ref{eq:asymtheta}), we get
\begin{equation}
F_n(c) = \big(1-\vartheta(\gamma^\star_n(c))\big)^n = \left(1-\frac{1+o(1)}{e\gamma^\star_n(c)}\right)^n.
\end{equation}
Now, using the
asymptotics in (\ref{eq:reftrancl}), we obtain that \whp
\begin{equation}
F_n = \begin{cases}
\big(1 - \frac{1+o(1)}{n^c \log n}\big)^n & \text{if $c < 1$} \\
\\[-1em]
\big(1 - \frac{1+o(1)}{2n \log \log n/\log n}\big)^n & \text{if
  $c=1$} \\
\\[-1em]
\big(1 - \frac{1+o(1)}{(1-c^{-1}) n}\big)^n & \text{if $c > 1$}
\end{cases}.\label{eq:phasetrans_nolimits}
\end{equation}
Upon taking the limit as $n \to \infty$ of the expressions given in
the right-hand side of (\ref{eq:phasetrans_nolimits}) and using
Lemma~\ref{lm:limits}, we obtain (\ref{eq:phasetrans}), and the
theorem is proved. 
\end{proof}

\subsection{Discussion}
\label{ssec:discuss}

In retrospect, it is easy to see that Theorem~\ref{thm:main} holds
trivially under the independent failure model. In order to appreciate
nontrivial features that appear once we pass to Pippenger's model, we
will invoke the game-theoretic interpretation given at
the end of Section~\ref{ssec:boole}.

Consider a Programmer-Hacker game of the kind described in
Section~\ref{ssec:boole}. Provided that the constant
$\varepsilon_0$ picked by the Referee is sufficiently small, the
Programmer has a good chance of winning if he sticks to the following
strategy: Let $n_0$ be the smallest integer and $c$ the largest
positive number, such that $\varepsilon_0 \le [e(1-1/c)n_0]^{-1}$ and $cn^{-1}_0\log n_0 \le 1$. The Programmer generates a random acyclic
digraph $\G_d(n_0,cn^{-1}_0 \log n_0)$ and uses it to construct a
Boolean network ${\sf N}_c(n)$ by adding variable inputs, assigning gate labels, and
designating the output gate, possibly in a completely arbitrary
fashion. He then hands this network to the Hacker. (Note that both the
Programmer and the Hacker have all the information needed to determine
which function $f$ is computed by the network.) If $c$ is large
enough, then Theorem~\ref{thm:main}
thus guarantees the existence of some $\delta < 1/2$
such that, with probability $1-o(1)$, the network generated by the
Programmer will $(\varepsilon_0, \delta)$-compute $f$, regardless of
the Hacker's actions. More precisely, for any $\varkappa,\varkappa' > 0$, there
exists large enough $N$ such that, for any $c \in (1,N/\log N)$ and
for all $n \ge N$, 
\begin{equation}
\Pr \left( 
\begin{array}{c}
\text{\small ${\sf N}_c(n)$ $([e (1-1/c)n]^{-1},\delta)$-computes $f$}\\
\text{\small with $\left|(1-\delta) - \exp\left(-\frac{c}{e (c-1)}\right)\right| <
  \varkappa$}
\end{array}
\right) >
1-\varkappa'.
\end{equation}
Provided that $e^{-c/[e(c-1)]} - \varkappa > 1/2$ and $\varepsilon_0
\le [e(1-1/c)n_0]^{-1}$ for some $n_0 \ge N$, the Programmer
will win with probability at least $1-\varkappa'$.

\section{Pippenger's model and stochastic domination}
\label{sec:stochdom}

In Section~\ref{ssec:motivation} an argument based on the
lopsided \Lovasz\ local lemma (\cite{ES91}, see also
Theorem~\ref{thm:4l} in this paper) allowed us to pinpoint the
possibility for a phase transition in Pippenger's model on random
graphs. In this section we show that the machinery of stochastic
domination \cite{Lig85,LSS97,Lin92} enables us to carry out a more refined
analysis of Pippenger's model. [We hasten to note that many of the
issues which we will touch upon have, in fact,
already been discussed by Pippenger in Ref.~\cite{Pip89}, but without
any systematic emphasis on stochastic domination.]

\subsection{Stochastic domination: the basics}
\label{ssec:dominate}
Once again, consider a finite acyclic digraph $G=(V,E)$ along
with the measurable space $(\Omega,\cF)$, where $\Omega = \B^V$ and
$\cF$ is the set of all subsets of $\Omega$. Elements of $\Omega$ are
binary strings of length $v(G)$; we will denote the $i$th component
(bit) of $\omega \in \Omega$ by $\omega(i)$. The total ordering of
$\B$ induces the following partial order of $\Omega$:
\begin{equation}
\omega_1 \prec \omega_2 \quad \Longleftrightarrow \quad \omega_1(i) \le
\omega_2(i), \forall i \in V.
\end{equation}
We say that a function $f : \Omega \to \R$ is {\em increasing} if
$\omega_1 \prec \omega_2$ implies $f(\omega_1) \le f(\omega_2)$. An
event $H \in \cF$ is called increasing if its indicator function,
$\idty_H$, is increasing. Informally speaking, an event is increasing
if its occurrence is unaffected by changing some bits from zero to
one. Decreasing functions and events are defined in an obvious way. Given two probability measures $\mu,\nu$ on $(\Omega,\cF)$, we
will say that $\mu$ is {\em stochastically dominated} by $\nu$ (and write
$\mu \preceq_s \nu$) if, for every increasing function $f$,
$\Ex_\mu(f) \le \Ex_\nu(f)$. As usual, the expectation $\Ex_\mu(\cdot)$ is
defined by
\begin{equation}
\Ex_\mu(f) := \int_\Omega f d\mu \equiv \sum_{\omega \in
  \Omega}f(\omega)\mu(\omega).
\end{equation}

Any probability measure $\mu$ on $(\Omega,\cF)$ is
equivalent to a family $\{ X_i : i \in V\}$ of $\B$-valued random
variables via
\begin{equation}
\Pr \big( \bigwedge_{i \in V}(X_i = \omega(i))\big) = \mu(\omega),
\qquad \forall \omega \in \Omega
\label{eq:jtlaw}
\end{equation}
(also cf. Section~\ref{ssec:pippenger}). We will say that the $X_i$
    {\em have joint law} $\mu$ if (\ref{eq:jtlaw}) holds. Then we have
    the following necessary and sufficient condition, due to Strassen
    \cite{Str65}, for one measure to dominate another (this is, in
    fact, a type of result that is proved most naturally by means of the
    so-called {\em coupling method}; see the monograph by Lindvall
    \cite{Lin92} for this as well as for many other
    useful applications of coupling).

\begin{lemma}Let $\mu,\nu$ be probability
    measures on $(\Omega,\cF)$. Then $\mu \preceq_s \nu$ if and only
    if there exist families of random variables $\{ X_i : i \in V\}$ and
    $\{ Y_i : i \in V\}$, defined on a common probability space, with
    respective joint laws $\mu$ and $\nu$, such that $X_i \le Y_i$
    almost surely for each $i \in V$.\label{lm:strassen}
\end{lemma}

Next we need a sufficient condition for a given probability
measure $\mu$ to dominate the Bernoulli product measure
$\pi^V_\eta$.  The following lemma is standard, and can be proved
along the lines of Holley \cite{Hol74} and Preston \cite{Pre74}:

\begin{lemma} Consider a family $\{ X_i : i \in V\}$ of random
    variables with joint law $\mu$. Suppose that there exists a total
    ordering $<$ of $V$ such that, for any $i \in V$ and any two disjoint
    sets $Y,Y' \subseteq V\backslash\{i\}$ with $j < i$ for all $j \in
    Y \cup Y'$, we have
\begin{equation}
\Pr_\mu \big( X_i = 1 | \bigwedge_{j \in Y} (X_j = 1) \wedge
\bigwedge_{j \in Y'}(X_j = 0) \big) \ge \eta,
\label{eq:holley}
\end{equation}
whenever $\Pr_\mu \big(\bigwedge_{j \in Y} (X_j = 1) \wedge
\bigwedge_{j \in Y'}(X_j = 0) \big) > 0$. Then $\mu \succeq_s
\pi^V_\eta$.
\label{lm:holley}
\end{lemma}
Note that Lemma~\ref{lm:holley} can also go in the other direction:
that is, if instead we write the $\le$ sign in (\ref{eq:holley}), then
we will have $\mu \preceq_s \pi^V_\eta$. 

Before we go on, let us introduce one more definition. Given $\eta \in
  (0,1)$, let us define a class $\cW^G_\eta$ of probability measures on
  $(\Omega,\cF)$ as follows. Let $\{ X_i : i\in V\}$ be the family of
  $\B$-valued random variables with joint law $\mu$. Then $\mu \in
  \cW^G_\eta$ if, for any $i \in V$ and for any disjoint sets $Y,Y'
  \subseteq V\backslash \bar{N}(i)$,
\begin{equation}
\mu \big(X_i = 1 | \bigwedge_{j \in Y}(X_j = 1) \wedge \bigwedge_{j \in
  Y'}(X_j = 0)\big) \ge \eta,
\label{eq:dubya}
\end{equation}
whenever the event we condition on has positive probability.

\subsection{Stochastic domination in Pippenger's model}
\label{ssec:pipstochdom}

In this section we will use the machinery of stochastic domination to
present a more delicate analysis of the $\varepsilon$-admissible model
of Pippenger. Let us fix (yet another!) acyclic digraph
$G=(V,E)$, which we will think of as a gate interconnection graph of a suitable
class of Boolean networks equivalent modulo gate labels and the
location of the output vertex (see Section~\ref{ssec:motivation} for
details). Given a probability measure $\mu \in \cM^G_\varepsilon$, let
$\{X_i : i \in V\}$ be a family of $\B$-valued random variables with
joint law $\mu$. Let us define $\tilde{X}_i := 1-X_i$ for each $i \in
V$; clearly, $\tilde{X}_i$ is an indicator random variable for the event
$\overline{A_i} \equiv \{ \text{gate $i$ does not fail} \}$. Let $\eta
:=1-\varepsilon$. The joint
law $\tilde{\mu}$ of $\{\tilde{X}_i\}$ is such that, for any $i \in V$
and any two disjoint sets $Y,Y' \in V\backslash \Gamma^\star(i)$,
\begin{equation}
\tilde{\mu} \big(\tilde{X}_i = 1 | \bigwedge_{j \in Y}(\tilde{X}_j = 1) \wedge
\bigwedge_{j\in Y'}(\tilde{X}_j=0)\big) \ge \eta,
\label{eq:notpipp}
\end{equation}
whenever the event we condition on has nonzero probability. Passing to the transitive
closure graph $G^\star = (V,E^\star)$, we can rewrite
(\ref{eq:notpipp}) as follows: for any $i \in V$ and any
disjoint $Y,Y'\subseteq V \backslash \bar{N}^\star(i)$,
\begin{equation}
\tilde{\mu} \big(\tilde{X}_i = 1 | \bigwedge_{j \in Y}(\tilde{X}_j = 1) \wedge
\bigwedge_{j \in Y'}(\tilde{X}_j = 0)\big) \ge \eta.
\label{eq:notpipp2}
\end{equation}
Comparing (\ref{eq:notpipp2}) and (\ref{eq:dubya}), we see that $\mu
\in \cM^G_\varepsilon$ implies $\tilde{\mu} \in
\cW^{G^\star}_\eta$. 

Let us show now that, for any $\mu \in \cM^G_\varepsilon$, the
corresponding $\tilde{\mu}$ satisfies $\tilde{\mu} \succeq_s
\pi^V_\eta$. We can use the same strategy as in the proof of
Lemma~\ref{lm:f_exact}. Namely, if we rearrange the vertices of $G$
according to some linear extension of the partial order $\preceq$,
then $\tilde{\mu}$ is easily seen to satisfy the conditions of
Lemma~\ref{lm:holley}, and we obtain the claimed result. It follows
directly from definitions that we have also $\mu \preceq_s
\pi^V_\varepsilon$. We also point out that Strassen's theorem
\cite{Str65} (Lemma~\ref{lm:strassen} in this paper) can be used to
give an amusing interpretation of Pippenger's model in terms of an
intelligent agent (``demon'') who, when faced with a realization of
$\{X_i\}$ with joint law $\mu \in \cM^G_\varepsilon$, can transform it
into a realization of random variables i.i.d. according to
$\pi^V_\varepsilon$ by changing some bits from 0 to 1, but
none from 1 to 0. (The same observation has been first made by
Pippenger \cite{Pip89}, who substantiated it using a
non-probabilistic result of Hwang \cite{Hwa79}.)

As part of our proof of Theorem~\ref{thm:main}, we have obtained the
exact formula $F_\varepsilon(G) = (1-\varepsilon)^{v(G)}$ by arranging
the vertices of $G$ according to a linear extension of the partial
order $\preceq$. The same conclusion can be easily reached once we
have established that, for any $\mu \in \cM^G_\varepsilon$, we have $\mu \preceq_s \pi^V_\varepsilon$ [or, equivalently, that the
  corresponding $\tilde{\mu} \in \cW^{G^\star}_\eta$ stochastically
  dominates $\pi^V_\eta$]; this extra
information enables us to obtain many other useful estimates
besides the one for $F_\varepsilon(G)$.

It is easy to see, for instance, that most of the ``really
interesting'' events one would naturally associate with proper operation of the
network (e.g., the event that the majority of the gates have not failed)
are decreasing events. That is, if such an event occurs in a
particular configuration $\omega$, then this event can be
destroyed by introducing additional failed gates. The event that no
gate fails is a particularly drastic example: it is destroyed if we
change the status of even a single gate. Equivalently, we may pass to
the corresponding probability measure $\tilde{\mu}$. In that case, given a
configuration $\omega \in \B^V$, the gate failures will correspond to
{\em zero} bits of $\omega$, with the nonzero bits indicating the gates that have not
failed. Therefore we may consider {\em increasing} events
if we agree to work with $\tilde{\mu}$ instead of $\mu$. Using the
stochastic domination properties of $\tilde{\mu}$, we get that
$\tilde{\mu}(H) \ge \pi^V_\eta(H)$ for any increasing event $H$.

As a simple example, consider the following set-up. Suppose we are given a network ${\sf
  N}$ whose output is the output of a gate that computes a Boolean function
$\varphi : \B^d \to \B$. Suppose that the inputs to this gate come
  from the outputs of subnetworks ${\sf N}_1,\ldots,{\sf N}_d$ (note
  that these subnetworks may, in general, share both gates and
  wires). Let $M$ be the event that the majority of the gates in ${\sf N}$
  have not failed, let $M_i$, $1\le i \le d$, be the event that the
  majority of the gates in ${\sf N}_i$ have not failed, and let $L$ be
  the event that the output gate of ${\sf N}$ has not failed. Then $L
  \cap \bigcap^d_{i=1}M_i$ implies $M$, so that
\begin{equation}
\Pr(M) \ge \Pr \big(L
  \cap \bigcap\nolimits^d_{i=1}M_i\big) = \Pr\big(L |
  \bigcap\nolimits^d_{i=1}M_i \big) \Pr \big(
  \bigcap\nolimits^d_{i=1}M_i \big).
\label{eq:pmaj}
\end{equation}
Suppose that the underlying probability measure $\mu$ is
$\varepsilon$-admissible. Then
\begin{equation}
\Pr\big(L | \bigcap\nolimits^d_{i=1}M_i \big) \ge \eta.
\label{eq:pmaj2}
\end{equation}
Likewise by $\varepsilon$-admissibility, $\tilde{\mu} \succeq_s
\pi^{\sf N}_\eta$ (see Section~\ref{ssec:boole} for this
notation). Therefore, since an
intersection of increasing events is increasing, we have
\begin{equation}
\Pr_{\tilde{\mu}}\big(
  \bigcap\nolimits^d_{i=1}M_i \big) \ge \pi^{\sf N}_\eta \big(
  \bigcap\nolimits^d_{i=1}M_i \big).
\label{eq:pmaj3}
\end{equation}
The right-hand side of (\ref{eq:pmaj3}) can be bounded from below by
means of the FKG inequality \cite{FKG71} (proved earlier by Harris
\cite{Har60} in the context of percolation) to give
\begin{equation}
\Pr_{\tilde{\mu}}\big(
  \bigcap\nolimits^d_{i=1}M_i \big) \ge \prod^d_{i=1} \pi^{{\sf N}_i}_\eta
  (M_i).
\label{eq:pmaj4}
\end{equation}
Let $a$ be the number of gates in the smallest of the subnetworks
${\sf N}_1,\ldots,{\sf N}_d$. Then, assuming that $\varepsilon \le
1/2$, Azuma-Hoeffding inequality \cite{Azu67,Hoe63} gives
\begin{equation}
\pi^{{\sf N}_i}_\eta(M_i) \ge 1-e^{-a (4\eta - 2)^2}.
\end{equation}
Putting everything together, we get
\begin{equation}
\Pr(M) \ge \eta \left(1-e^{-a (4\eta-2)^2} \right)^d.
\end{equation}
This use of the FKG inequality is similar to that in Feder
\cite{Pip89}, whose work was concerned with the depth and reliability
bounds for reliable computation of Boolean functions under the
independent failure model.

At this point we also mention another PFM discussed by Pippenger in
Ref.~\cite{Pip89} --- namely, the so-called $\varepsilon$-{\em majorized}
model. Under this model, a network ${\sf N}$ is mapped to the set of
all probability measures $\mu$ on $\B^{\hat{V}_{\sf N}}$ that are
stochastically dominated by the Bernoulli product measure $\pi^{\sf
  N}_\varepsilon$. It follows easily from the discussion above that
the $\varepsilon$-majorized model is more stringent than the
$\varepsilon$-admissible model. As an example of its use, we can
mention the work of Dobrushin and Ortyukov \cite{DO77b}, where a
result proven for the $\varepsilon$-majorized model automatically
carries over to the $\varepsilon$-independent one.

\section{Closing remarks and future directions}
\label{sec:close}

In this paper we have showed that a phase transition is possible in
the $\varepsilon$-admissible PFM of Pippenger \cite{Pip89} as soon as
random graphs show up in the picture. In doing so, we have barely
scratched the surface of a wonderfully rich subject --- namely, the
statistical mechanics of multicomponent random systems on directed
graphs. Most of the work connected to phase transitions in large
combinatorial structures has been done in the context of undirected
graphs, since the methods of
statistical mechanics applicable to the study of combinatorial
structures have been originally developed in that context as well.

For instance, the independent-set polynomial of a simple
undirected graph \cite{FS90,HL94,She85} can also be viewed as a partition
function of a repulsive lattice gas \cite{Sim93}, and the powerful
machinery of cluster expansions \cite{Dob96} developed in the latter
setting can also be applied quite successfully to the former; the reader is
encouraged to consult a recent paper by Scott and Sokal \cite{SS03}
for an exposition of these matters from the
viewpoint of both statistical mechanics and graph theory. In the
future we would like to study applications
of statistical mechanics to combinatorial structures with
directionality. So far, very few results along these lines have
appeared; the papers of
Whittle \cite{Whi89,Whi90} are among the few examples known to the
present author where a partition function is constructed for a class
of statistical-mechanical models on directed graphs. The relative
dearth of applications of statistical mechanics to structures with
directionality is due to the fact that, once directionality is
introduced, the symmetry needed for the lattice-gas formalism is
destroyed, and it is not immediately evident how one could relate
combinatorial properties of directed graphs to mathematical objects of
statistical mechanics. (Incidentally, this very point has also been
brought up by Scott and Sokal \cite{SS03}.)

\renewcommand{\theequation}{A.\arabic{equation}}
\renewcommand{\thetheorem}{A.\arabic{theorem}}
\setcounter{equation}{0}
\setcounter{theorem}{0}
\section*{Appendix}

Our goal in this appendix is to prove a theorem that can be thought of as
a specialization of the lopsided \Lovasz\ local lemma of \Erdos\ and
Spencer \cite{ES91} (see also Theorem~\ref{thm:4l}) to families of
random variables whose joint laws are elements of $\cW^G_\eta$. Both
the theorem and its proof go very much along the lines a similar
result of Liggett, Schonmann, and Stacey \cite{LSS97}, except that
theirs was formulated for undirected graphs. 

\begin{theorem}\label{thm:lllsupercharge} Let $G=(V,E)$ be a directed
  acyclic graph, in which every vertex has out-degree at most $\Delta
  \ge 1$. Let $\mu \in \cW^G_\eta$ with $\varepsilon := 1-\eta \le
  \Delta^\Delta/(\Delta+1)^{(\Delta+1)}$. Then $\mu \succeq_s
  \pi^V_\rho$, where
\begin{equation}
\rho = \left( 1-
\frac{\varepsilon^{1/(\Delta+1)}}{\Delta^{\Delta/(\Delta+1)}}
\right)\left(1-(\varepsilon \Delta)^{1/(\Delta+1)}\right).
\label{eq:rho}
\end{equation}
\end{theorem}

\begin{proof}First we need a lemma.
\begin{lemma}\label{lm:xyz}Let $G=(V,E)$ satisfy the conditions of
  Theorem~\ref{thm:lllsupercharge}. Given $\eta \in (0,1)$, consider
  $\mu \in \cW^G_\eta$. Suppose that there exist constants
  $\alpha,\lambda \in (0,1)$, such that
\begin{equation}
\varepsilon \le (1-\alpha)(1-\lambda)^\Delta,
\label{eq:cond1}
\end{equation}
\begin{equation}
\varepsilon \le (1-\alpha)\alpha^\Delta.
\label{eq:cond2}
\end{equation}
Consider a family $\{X_i : i \in V\}$ of random variables with joint
law $\mu$, and let $\{ Y_i : i \in V\}$ be a family of random
variables, independent of $\{ X_i\}$ and with joint law
$\pi^V_\lambda$. Let $Z_i := X_iY_i$ for each $i \in V$. Then, for
each $i \in V$, each $Y \subseteq V\backslash\{i\}$, and each $z \in
\B^{|Y|}$, we have
\begin{equation}
\Pr \big(Z_i = 1 | \bigwedge\nolimits^{|Y|}_{j=1}(Z_{i_j} = z_j) \big) \ge
\alpha \lambda,
\label{eq:zdom}
\end{equation}
where the $i_j$ are elements of $Y$.
\end{lemma}

\begin{remark} The corresponding theorem of Liggett,
  Schonmann, and Stacey \cite{LSS97} is formulated in terms of an {\em
  undirected} graph $G$, with $\Delta$ being the maximum {\em degree}
  of a vertex. Therefore they need to impose an additional condition,
  namely that $i$ is adjacent to at most $\Delta-1$ vertices in $Y$.
  As a consequence, one has to make the replacement $\Delta \to \Delta
  - 1$, e.g., in (\ref{eq:rho}), (\ref{eq:cond1}) and
  (\ref{eq:cond2}). However,
  because here we deal with {\em directed} graphs and $\Delta$ is the
  maximum {\em out-degree} of a vertex, there are automatically no more than
  $\Delta$ vertices $j$ in $G$ with $(i,j) \in E$.
\end{remark}

\begin{proof} Note that (\ref{eq:zdom}) is equivalent to
\begin{equation}
\Pr \big(X_i = 1 | \bigwedge\nolimits^{|Y|}_{j=1}(Z_{i_j} = z_j) \big) \ge \alpha
\label{eq:zdom_equiv}
\end{equation}
due to independence of $\{X_i\}$ and $\{Y_i\}$ and to the fact that
$\lambda > 0$. We will proceed by proving (\ref{eq:zdom_equiv}) by
induction on $|Y|$.

Suppose first that $Y=\varnothing$. Then (\ref{eq:zdom_equiv}) is simply the
statement that $\Pr (X_i=1) \ge \alpha$. Now, $\Pr (X_i=1) \ge \eta$
because $\mu \in \cW^G_\eta$, and $\eta \ge \alpha$ by
(\ref{eq:cond2}). Thus suppose that (\ref{eq:zdom_equiv}) holds for
all $Y \subseteq V\backslash\{i\}$ with $|Y| < J$,
where $J \ge 1$. We will prove that it also holds for $|Y| = J$.

Fix $Y = \{ i_1,\ldots,i_J \}$ and $z \in \B^J$. We write $Y$ as a
disjoint union $M_0 \cup M_1 \cup M$, where
\begin{eqnarray}
M_0 &:=& \{ i_j, 1 \le j \le J : \text{$(i,i_j) \in E$ and $z_j =
  0$}\} \\
M_1 &:=& \{ i_j, 1 \le j \le J : \text{$(i,i_j) \in E$ and $z_j =
  1$}\} \\
M &:=& Y \backslash (M_0 \cup M_1).
\end{eqnarray}
Let us also define the events
\begin{eqnarray}
A_0 &:=& \{ Z_{i_j} = 0 : i_j \in M_0 \} \\
B_0 &:=& \{ Y_{i_j} = 0 : i_j \in M_0 \} \\
A_1 &:=& \{ X_{i_j} = 1 : i_j \in M_1 \} \\
A &:=& \{ Z_{i_j} = z_j : i_j \in N \} 
\end{eqnarray}
Now, for any $j \in V$, $Z_j = 1$ is by definition equivalent to both
$X_j = 1$ and $Y_j = 1$. Furthermore, $\{ X_j \}$ and $\{Y_j\}$ are
independent. Therefore we can write
\begin{equation}
\Pr \big(X_i = 1 | \bigwedge\nolimits^J_{j=1}(Z_{i_j} = z_j) \big)
= \Pr \big(X_i = 1 | A_0 \cap A_1 \cap A \big)
\end{equation}
Now
\begin{eqnarray}
 \Pr \big(X_i = 1 | A_0 \cap A_1 \cap A \big) &=& 1 -  \Pr \big(X_i =
 0 | A_0 \cap A_1 \cap A \big) \nonumber \\
&=& 1 - \frac{\Pr\big(X_i = 0, A_0 \cap A_1 \cap A\big)}{\Pr\big(A_0
 \cap A_1 \cap A\big)} \nonumber \\
&\ge& 1- \frac{\Pr\big(X_i=0, A\big)}{\Pr\big(B_0 \cap A_1 \cap
 A\big)} \nonumber \\
&= & 1 - \frac{\Pr\big(X_i=0 | A\big)}{\Pr\big(B_0 \cap A_1 | A\big)}.\label{eq:cond_bound}
\end{eqnarray}
Since $(i,i_j) \not\in E$ for all $i_j \in M$, the numerator is at
most $\varepsilon$. The denominator is equal to $(1-\lambda)^{|M_0|}
\Pr\big(A_1 | A\big)$. Assume that $M_1 = \{k_1,\ldots,k_s\}$, $s =
|M_1|$. Then,
\begin{eqnarray}
\Pr\big(A_1 | A\big) &=& \Pr\big(\bigwedge\nolimits^s_{\ell =
  1}(X_{k_\ell} = 1) | A\big) \nonumber\\
&=& \Pr\big(X_{k_1} = 1 | A\big) \prod^{s-1}_{\ell = 1} \Pr
  \big(X_{k_{\ell + 1}} = 1 | A, \bigwedge\nolimits^\ell_{m=1} (X_{k_m}
  = 1) \big) \nonumber\\
&\ge& \alpha^{|M_1|},
\end{eqnarray}
where in the last step we have applied the inductive
hypothesis to each of the terms in the product. Therefore
\begin{equation}
 \Pr \big(X_i = 1 | A_0 \cap A_1 \cap A \big) \ge 1 -
 \frac{\varepsilon}{(1-\lambda)^{|M_0|}\alpha^{|M_1|}}.
\label{eq:almost_there}
\end{equation}
Since $|M_0|+|M_1| \le \Delta$ by hypothesis, and $\varepsilon/(1-\alpha) \le
\alpha^\Delta \le 1$ by (\ref{eq:cond2}), we have
\begin{equation}
(1-\lambda)^{|M_0|}\alpha^{|M_1|} \ge
  \left(\frac{\varepsilon}{1-\alpha}\right)^{(|M_0|+|M_1|)/\Delta} \ge \frac{\varepsilon}{1-\alpha}.
\end{equation}
Therefore the right-hand side of
(\ref{eq:almost_there}) is at least $1-\varepsilon/[\varepsilon/(1-\alpha)] =
\alpha$, and the lemma is proved.
\end{proof}

Now let $\{X_i\}$, $\{Y_i\}$, and $\{Z_i\}$ be as in
Lemma~\ref{lm:xyz}. Let $\nu$ be the joint law of $\{Z_i\}$. By
construction, $Z_i \le X_i Y_i$ for each $i \in V$, so $\mu \succeq_s
\nu$ by Lemma~\ref{lm:strassen}. We now show that $\nu \succeq_s
\pi^V_{\alpha \lambda}$. Let $<$ be an {\em arbitrary} total ordering
of $V$. Then Lemma~\ref{lm:xyz} and Lemma~\ref{lm:holley} imply that
$\nu \succeq_s \pi^V_{\alpha \lambda}$. Thus $\mu \succeq_s
\pi^V_{\alpha \lambda}$.

To conclude the proof, suppose that
\begin{equation}
\varepsilon \le \frac{
\Delta^\Delta}{(\Delta+1)^{\Delta + 1}}.
\label{eq:epsicond}
\end{equation}
Let
\begin{equation}
\alpha = 1- \frac{\varepsilon^{1/(\Delta+1)}}{\Delta^{\Delta/(\Delta+1)}}
\quad \text{and} \quad \lambda = 1-(\varepsilon
\Delta)^{1/(\Delta+1)}.
\label{eq:alpha_lambda}
\end{equation}
Then $(1-\alpha)(1-\lambda)^\Delta = \varepsilon$, which yields
(\ref{eq:cond1}). Condition (\ref{eq:epsicond}) is equivalent to
\begin{equation}
\varepsilon^{1/(\Delta+1)} \le
\frac{\Delta^{\Delta/(\Delta+1)}}{\Delta+1},
\end{equation}
which, when substitued into (\ref{eq:alpha_lambda}), yields
\begin{equation}
\alpha \ge 1 - \frac{1}{\Delta+1} \quad \text{and} \quad \lambda \ge
\frac{1}{\Delta+1}.
\end{equation}
This shows that the choice we have made in (\ref{eq:alpha_lambda})
leads to $\alpha,\lambda \in [0,1]$, and that $1-\lambda \ge
\alpha$. The latter inequality implies (\ref{eq:cond2}). Therefore, by
Lemma~\ref{lm:xyz}, $\mu \succeq_s \pi^V_{\alpha\lambda}$, and the
theorem is proved.
\end{proof}

It is important to mention that Theorem~\ref{thm:lllsupercharge} is
useful only when $G$ is not transitively closed, i.e., when $(i,j) \in
E$ and $(j,k) \in E$ does not imply $(i,k) \in E$. Otherwise one can
partially order the vertices of $G$ by defining $i \preccurlyeq j$ if
$(i,j) \in E$ for distinct $i,j \in V$, and $i \preccurlyeq i$ for each $i
\in V$. As usual, let $i_1 < \ldots < i_{v(G)}$ be a total order of $V$
according to some linear extension of $\preccurlyeq$. Thus, for any
$i_j$, all the $i_k$ with $i_k < i_j$ and distinct from $i_j$ are not
in $\bar{N}(i_j)$. Therefore we can apply Lemma~\ref{lm:holley}
directly to any $\mu \in \cW^G_\eta$ to conclude that $\mu \succeq_s
\pi^V_\eta$. This is, in fact, precisely the case we have dealt with
in this paper --- namely, when $G$ is a transitive closure of some
other acyclic digraph $G_0$.

\section*{Acknowledgments} I would like to thank Svetlana Lazebnik for
helpful discussions.


\begin{thebibliography}{99}
\small

\bibitem{Aig79}M.~Aigner, \emph{Combinatorial Theory},
  Springer-Verlag, Berlin, 1979. 

\bibitem{AB02}R.~Albert and A.-L.~Barab\'asi, Statistical mechanics of
  complex networks, \emph{Rev. Mod. Phys.} {\bf 74} (2002), 47--97. 

\bibitem{AS00}N.~Alon and J.H.~Spencer, \emph{The Probabilistic
  Method}, 2nd ed., Wiley, New York, 2000. 

\bibitem{Azu67}K.~Azuma, Weighted sums of certain dependent variables,
\emph{T\^ohoku Math. J.} {\bf 3} (1967), pp. 357--367. 

\bibitem{BE84}A.~Barak and P.~\Erdos, On the maximal number of
  strongly independent vertices in a random acyclic directed graph,
  \emph{SIAM J. Algebraic and Discrete Methods} {\bf 5} (1984),
  508--514. 

\bibitem{BMW00}G.~Biroli, R.~Monasson, and M.~Weigt, A variational
  description of the ground state structure in random satisfiability
  problems, \emph{Eur. Phys. J. B} {\bf 14} (2000), 551--568. 

\bibitem{Bol01}B.~Bollob\'as, \emph{Random Graphs}, 2nd ed., Cambridge
  University Press, Cambridge, 2001. 

\bibitem{BBCKW01}B.~Bollob\'as, C.~Borgs, J.T.~Chayes, J.H.~Kim, and
  D.B.~Wilson, The scaling window of the 2-SAT transition, \emph{Random
  Struct. Alg.} {\bf 18} (2001), 201--256. 

\bibitem{BCP01}C.~Borgs, J.T.~Chayes, and B.~Pittel, Phase transition
  and finite-size scaling for the integer partitioning problem,
  \emph{Random Struct. Alg.} {\bf 19} (2001), 247--288.

\bibitem{Dob96}R.L.~Dobrushin, Estimates of semi-invariants for the
  Ising model at low temperatures, \emph{Topics in Statistical and Theoretical Physics}, American Mathematical
Society Translations, Ser. 2, vol.~177 (1996), 59--81.

\bibitem{DO77a}R.L.~Dobrushin and S.I.~Ortyukov, Lower bound for the
  redundancy of self-correcting arrangements of unreliable functional
  elements, \emph{Prob. Inf. Transm.} {\bf 13} (1977), 59--65. 

\bibitem{DO77b}R.L.~Dobrushin and S.I.~Ortyukov, Upper bound on the
redundancy of self-correcting arrangements of unreliable functional
elements, \emph{Prob. Inf. Transm.} {\bf 13} (1977), 203--218. 

\bibitem{Dur96}R.~Durrett, \emph{Probability: Theory and Examples},
  2nd ed., Wadsworth, Belmont, 1996. 

\bibitem{EL75}P.~\Erdos\ and L.~\Lovasz, Problems and results on
  3-chromatic hypergraphs and some related questions, \emph{Infinite
  and Finite Sets}, A.~Hajnal \emph{et al.}, eds., North-Holland
  (1975), 609--628. 

\bibitem{ES91}P.~\Erdos\ and J.~Spencer, Lopsided \Lovasz\ local lemma
  and Latin transversals, \emph{Discrete Appl. Math.} {\bf 30} (1991),
  151--154. 

\bibitem{Fed89}T.~Feder, Reliable computation by networks in the
  presence of noise, \emph{IEEE Trans. Inform. Theory} {\bf 35}
  (1989), 569--571. 

\bibitem{FS90}D.C.~Fisher and A.E.~Solow, Dependence polynomials,
  \emph{Discrete Math.} {\bf 82} (1990), 251--258.

\bibitem{FKG71}C.M.~Fortuin, P.W.~Kasteleyn, and J~Ginibre,
  Correlation inequalities on some partially ordered sets,
  \emph{Commun. Math. Phys.} {\bf 22} (1971), 89--103. 

\bibitem{GG94}P.~G\'acs and A.~G\'al, Lower bounds for the complexity
  of reliable Boolean circuits with noisy gates, \emph{IEEE
  Trans. Inform. Theory} {\bf 40} (1994), 579--583. 

\bibitem{Har60}T.E.~Harris, A lower bound on the critical probability
  in a certain percolation process, \emph{Proc. Cambridge Phil. Soc.}
  {\bf 56} (1960), 13--20. 

\bibitem{HL94}C.~Hoede and X.-L.~Li, Clique polynomials and
  independent set polynomials of graphs, \emph{Discrete Math.} {\bf
  125} (1994), 219--228.

\bibitem{Hoe63}W.~Hoeffding, Probability inequalities for sums of
bounded random variables, \emph{J. Amer. Stat. Assoc.} {\bf 58}
(1963), 13--30. 

\bibitem{Hol74}R.~Holley, Remarks on the FKG inequalities,
  \emph{Commun. Math. Phys.} {\bf 36} (1974), 227--231. 

\bibitem{Hwa79}F.K.~Hwang, Majorization on a partially ordered set,
  \emph{Proc. Amer. Math. Soc.} {\bf 76} (1979), 199--203. 

\bibitem{JLR00}S.~Janson, T.~\L uczak, and A.~Ruci\'nski, \emph{Random
  Graphs}, Wiley, New York, 2000. 

\bibitem{Lig85}T.M.~Liggett, \emph{Interacting Particle Systems},
  Springer, New York, 1985. 

\bibitem{LSS97}T.M.~Liggett, R.H.~Schonmann, and A.M.~Stacey,
  Domination by product measures, \emph{Ann. Probab.} {\bf 25} (1997),
  71--95. 

\bibitem{Lin92}T.~Lindvall, \emph{Lectures on the Coupling Method},
  Wiley, New York, 1992. 

\bibitem{LS97}B.~Luque and R.V.~Sol\'e, Phase transitions in random
  networks: simple analytic determination of critical points,
  \emph{Phys. Rev. E} {\bf 55} (1997), 257--260. 

\bibitem{MMZ01}O.C.~Martin, R.~Monasson, and R.~Zecchina, Statistical
  mechanics methods and phase transitions in optimization problems,
  \emph{Theoret. Comput. Sci.} {\bf 265} (2001), 3--67. 

\bibitem{Mul56}D.E.~Muller, Complexity in electronic switching
  circuits, \emph{IRE Trans. Electr. Comput.} {\bf 5} (1956), 15--19. 

\bibitem{Neu56}J.~von Neumann, Probabilistic logics and the synthesis
  of reliable organisms from unreliable components, \emph{Automata
  Studies}, C.E. Shannon and J. McCarthy, eds., Princeton University
  Press (1956), 329--378. 

\bibitem{Pip88}N.~Pippenger, Reliable computation by formulas in the
  presence of noise, \emph{IEEE Trans. Inform. Theory} {\bf 34}
  (1988), 194--197. 

\bibitem{Pip89}N.~Pippenger, Invariance of complexity measures for
  networks with unreliable gates, \emph{J. Assoc. Comput. Machinery}
  {\bf 36} (1989), 531--539. 

\bibitem{Pip90}N.~Pippenger, Developments in ``The Synthesis of
  Reliable Organisms from Unreliable Components,'' \emph{Legacy of
  J.~von Neumann}, Proceedings of Symposia in Pure Mathematics,
  vol.~50 (1990), 311--324.

\bibitem{PST91}N.~Pippenger, G.D.~Stamoulis, and J.N.~Tsitsiklis, On a
  lower bound for the redundancy of reliable networks with noisy
  gates, \emph{IEEE Trans. Inform. Theory} {\bf 37} (1991), 639--643. 

\bibitem{PT01}B.~Pittel and R.~Tungol, A phase transition phenomenon
  in a random directed acyclic graph, \emph{Random Struct. Alg.} {\bf
  13} (2001), 164--184. 

\bibitem{Pre74}C.J~Prestion, A generalization of the FKG inequalities,
  \emph{Commun. Math. Phys.} {\bf 36} (1974), 233--242. 

\bibitem{SS03}A.D.~Scott and A.D.~Sokal, The repulsive lattice gas,
  the independent-set polynomial, and the \Lovasz\ local lemma,
  cond-mat/0309352 at arXiv.org. 

\bibitem{SC01}G.~Semerjian and L.F.~Cugliandolo, Cluster expansions in
  dilute systems: applications to satisfiability problems and spin
  glasses, \emph{Phys. Rev. E} {\bf 64} (2001), 036115. 

\bibitem{She85}J.B.~Shearer, On a problem of Spencer,
  \emph{Combinatorica} {\bf 5} (1985), 241--245. 

\bibitem{Sim93}B.~Simon, \emph{The Statistical Mechanics of Lattice
  Gases}, Princeton University Press, Princeton, 1993.

\bibitem{Str65}V.~Strassen, The existence of probability measures with
  given marginals, \emph{Ann. Math. Statist.} {\bf 36} (1965),
  423--439. 

\bibitem{Weg87}I.~Wegener, \emph{The Complexity of Boolean Functions},
  Wiley, New York, 1987. 

\bibitem{Whi89}P.~Whittle, The statistics of random directed graphs,
  \emph{J. Stat. Phys.} {\bf 56} (1989), 499--516. 

\bibitem{Whi90}P.~Whittle, Fields and flows on random graphs,
  \emph{Disorder in physical systems}, G.R.~Grimmett and D.J.A.~Welsh,
  eds., Oxford University Press (1990), 337--348.

\end{thebibliography}
\end{document}